\documentclass{amsart}
\usepackage{amsfonts}
\newtheorem{theorem}{Theorem}

\newtheorem{rem}{Remark}
\newtheorem{rems}{Remarks}
\newtheorem{lm}{Lemma}
\newtheorem{ex}{Example}
\newtheorem{cor}{Corollary}
\newtheorem{prop}{Proposition}
\newtheorem{nota}{Notation}

\begin{document}

\title{A mapping defined by the Schur-Szeg\H{o} composition}
\author{Vladimir P. Kostov}
\address{Universit\'e de Nice, Laboratoire de
Math\'ematiques, Parc Valrose, 06108 Nice, Cedex 2, France}
\email{kostov@math.unice.fr}
\thanks{Research supported by the French Foundation CNRS under Project 20682.}
\subjclass[2000]{Primary 12D10; Secondary 30D99}
\keywords{Schur-Szeg\H{o} composition; composition factor; hyperbolic polynomial}
\maketitle 

\begin{abstract}
Each degree $n+k$ polynomial of the form $(x+1)^k(x^n+c_1x^{n-1}+\cdots +c_n)$, $k\in \mathbb{N}$,
is representable as Schur-Szeg\H{o} composition of $n$ polynomials of the form $(x+1)^{n+k-1}(x+a_j)$. We study
properties of the affine mapping $\Phi _{n,k}$~:~$(c_1,\ldots ,c_n)$ $\mapsto$
$(\sigma _1, \ldots ,\sigma _n)$, where $\sigma _i$ are the elementary symmetric polynomials of the numbers $a_j$.
We study also properties of a similar mapping for functions of the form $e^xP$, where $P$ is a polynomial, 
$P(0)=1$, and we extend the Descartes rule to them.
\end{abstract}

The {\em Schur-Szeg\H{o} composition (SSC)} of two polynomials $A:=\sum _{j=0}^n{n\choose j}\alpha _jx^j$
and $B:=\sum _{j=0}^n{n\choose j}\beta _jx^j$ is defined by the formula
$A*B:=\sum _{j=0}^n{n\choose j}\alpha _j\beta _jx^j$. The SSC is commutative and associative.
The above formula can be generalized in a self-evident way to the case of composition of more than two polynomials.

Obviously, $(x+1)^n*A=A$ for any degree $n$ polynomial $A$; that is, in the case of the SSC 
the polynomial $(x+1)^n$ plays the role of unity.
If the polynomials $A$ and $B$ are considered as degree
$n+k$ ones, their first $k$ coefficients being equal to $0$, then the formula for $A*B$ will be a different one.
To avoid such an ambiguity we assume throughout this paper that the leading coefficient of at least one
of the composed polynomials is non-zero. See more about the SSC in \cite{Pr}
and \cite{RS}.

In this paper we study the affine mappings $\Phi _{n,k}$ (connected with the SSC and 
defined after the proof of Lemma~\ref{zero}) 
and their generalization $\Phi$ for entire functions (defined before Remarks~\ref{limit}). 
We also generaize the Descartes rule, see Theorem~\ref{Descartes}. 
The following formulae are proved in \cite{Ko1} ($S$ is a degree
$n-1$ polynomial):

\begin{equation}\label{differ}
(A*B)'=(1/n)(A'*B')~~~~~~~,~~~~~~~~~(xS*B)=(x/n)(S*B')~~.
\end{equation}

\begin{prop} \label{mult}
{\rm (Proposition 1.4 in \cite{KoSh}.)}
If the polynomials $A$ and $B$ have roots $x_A\neq 0$ and $x_B\neq 0$ of multiplicities
$m_A$ and $m_B$ respectively, where $m_A+m_B\geq n$, then $-x_Ax_B$ is a root of $A*B$
of multiplicity $m_A+m_B-n$.
\end{prop}

%The conditions $x_A\neq 0$, $x_B\neq 0$ are omitted in \cite{KoSh} 
%which is wrong.
The following proposition is used to define below the mappings
$\Phi _{n,k}$, $k\geq 1$:

\begin{prop}\label{Kprop}
Each polynomial
$P:=(x+1)^k(x^n+c_1x^{n-1}+\cdots +c_n)$ is representable as
SSC 

\begin{equation}\label{Kformula}
P=K_{n,k;a_1}*\cdots *K_{n,k;a_n}~~~{\it with}~~~K_{n,k;a_i}:=(x+1)^{n+k-1}(x+a_i)~,
\end{equation}
where the complex numbers $a_i$ are unique up to permutation.
\end{prop}
\newpage 

{\em Proof:}\\ 

For $k=1$ the proposition is announced in Remark~7 of \cite{Ko1} and is proved in \cite{AlKo}.
For $k>1$ it can be deduced from there as
follows: write $P$ in the form $(x+1)((x+1)^{k-1}(x^n+c_1x^{n-1}+\cdots +c_n))$. The
second factor is a polynomial of degree $n+k-1$ to which one can apply the result from \cite{AlKo}
with $n$ replaced by $n+k-1$. Hence $P$ is SSC of $n+k-1$
{\em composition factors} $K_{n,k;a_i}$. One can deduce from Proposition~\ref{mult} that
$k-1$ of these composition factors equal $K_{n,k;1}$ (because $-1$ is a $(k-1)$-fold root of the second factor) and hence can be skipped.~~~~~$\Box$

\begin{lm}\label{zero}
The coefficient of $x^s$ in $P$ is zero if and only if one of the numbers $a_i$ equals $-s/(n+k-s)$.
\end{lm}

This follows from the formula

\begin{equation}\label{Keq}
K_{n,k;a_i}=
\sum_{s=0}^{n+k}{n+k \choose s}\left( \frac{n+k-s}{n+k}\, a_i+\frac{s}{n+k}\right) x^s~.
\end{equation}
Indeed, the coefficient of $x^s$ 
in at least one polynomial $K_{n,k;a_i}$ must equal $0$.~~~~~$\Box$\\

With $c_i$ and $a_i$ as in Proposition~\ref{Kprop}, the mapping $\Phi _{n,k}$ is defined
like this:
$$\Phi _{n,k}~:~(c_1,\ldots ,c_n)\mapsto (\sigma _1,\ldots ,\sigma _n)~,~~
{\rm where}~~\sigma _j:=\sum _{1\leq i_1<\cdots <i_j\leq n}a_{i_1}\cdots a_{i_j}~.$$

The mapping $\Phi _{n,k}$ is affine. For $k=1$ this is proved in \cite{Ko2}. For any $k$ it follows from there
(the coefficients of the polynomial $P/(x+1)$ are affine functions of the variables $c_i$). Properties of
$\Phi _{n,1}$ are studied in \cite{Ko2}, \cite{Ko3}, \cite{Ko4} and \cite{KoShMa}. In this
paper we continue the study of paper \cite{Ko4} and extend it to the case of entire functions.

The SSC of the entire functions $f:=\sum _{j=0}^{\infty}\gamma _jx^j/j!$ and 
$g:=\sum _{j=0}^{\infty}\delta _jx^j/j!$ is defined by the formula $f*g=\sum _{j=0}^{\infty}\gamma _j\delta _jx^j/j!$. 
Set $P_m:=1+c_{1}x+\cdots +c_mx^m$, $\tilde{\sigma}_k:=\sum _{1\leq j_1<\cdots <j_k\leq m}1/(a_{i_1}\cdots a_{i_k})$. 
The following proposition allows to define an analog of the mappings $\Phi _{n,k}$: 

\begin{prop}\label{entire}
%{\rm (Theorem~3 in \cite{DiKo}.)} 
Each function $e^xP_m$, where $P_m$ is a degree $m$ polynomial such that $P_m(0)=1$,
is representable in the form 

\begin{equation}\label{composePhi}
e^xP_m=\kappa _{a_1}*\cdots *\kappa _{a_m}~,~~~~~~~~{\it where}~~~~~~~~\kappa _{a_j}=e^x(1+x/a_j)~.
\end{equation}
The numbers $a_j$ are unique up to permutation.
\end{prop}

Indeed, it is easy to show by induction on $m$ (the proof is left for the reader) 
that the SSC of $m$ composition factors $\kappa _{a_j}$ is of the form 
$(1+\sum _{i=1}^mb_ix^i)e^x$, where $b_i=\sum _{l=i}^m\zeta _{i,l}\tilde{\sigma}_l$, $\zeta _{i,l}\in \mathbb{N}$, 
$\zeta _{i,i}=1$. The mapping $(\tilde{\sigma}_1,\ldots ,\tilde{\sigma}_m)\mapsto (b_1,\ldots ,b_m)$ 
is linear upper-triangular and non-degenerate from 
where the proposition follows.~~~~~$\Box$

Define the mapping $\Phi$ as follows: 
$\Phi ~:~(c_1,\ldots ,c_m)\mapsto (\tilde{\sigma}_1,\ldots ,\tilde{\sigma}_m)$. 

\begin{rems}\label{limit}
{\rm 1) The mapping $\Phi$ is a limit of mappings $\Phi _{n,k}$ as $k\rightarrow \infty$: 
each polynomial $k^k(x/k+1)^k(x^n+c_1x^{n-1}+\cdots +c_n)$ can be represented as 
SSC of $n$ {\em composition factors} of the form 
$k^{n+k-1}(x/k+1)^{n+k-1}(x+a_i)$. The 
proof of this is completely analogous to the proof of Proposition~\ref{Kprop}. There 
remains to observe that $\lim _{k\rightarrow \infty}(x/k+1)^k=e^x$. To avoid the constant factors $k^k$ and 
$k^{n+k-1}$ which tend to infinity as $k\rightarrow \infty$, one can consider instead polynomials 
$(x/k+1)^k(c_0x^n+c_1x^{n-1}+\cdots +c_{n-1}x+1)$ and composition factors of the form 
$(x/k+1)^{n+k-1}(x/a_i+1)$ or $(x/k+1)^{n+k-1}x$ in which case no constant factors are necessary. 

2) For the {\em composition factors} $\kappa _{a_i}$ one has the formula}

\begin{equation}\label{kappa}
\kappa _{a_j}=\sum _{j=0}^{\infty}(1/j!)(1+j/a_j)x^j~.
\end{equation}

{\rm 3) If $P$ (resp. $P_m$) is a real polynomial, then part of the numbers $a_j$ in formula (\ref{Kformula}) 
(resp. (\ref{composePhi})) are real and the rest form complex conjugate couples. Indeed, otherwise 
conjugation of the two sides of (\ref{Kformula}) or (\ref{composePhi}) would produce a new set of numbers 
$a_j$ which contradicts their uniqueness.}
\end{rems}

\begin{nota}
{\rm We denote by $U_n\subset \mathbb{R}^n\cong Oc_1\cdots c_n$  
the subset defined by the conditions $(-1)^ic_i\geq 0$. By $\Pi _{n}$ we denote 
the {\em hyperbolicity domain} of the family of polynomials $P$, i.e. the set of values of 
the coefficients $c_i$ for which $P$ is hyperbolic. We write $V_{n}\subset \mathbb{R}^{n}$ for 
the set of values of the coefficients of $P$ for which the real parts of all roots are non-negative. 
It is easy to show that $(\Pi _{n}\cap U_{n})\subset V_{n}\subset U_{n}$. By $T[f]$ we denote 
the Taylor series at $0$ of the entire function $f$.}
\end{nota}

\begin{theorem}\label{Utm}
For each $n\geq 1$ and for each $k\geq 1$ one has $\Phi _{n,k}(U_{n})\subset U_{n}$.
\end{theorem}

\begin{cor}\label{Ucor}
For the mapping $\Phi$ one has $\Phi (U_{n})\subset U_{n}$.
\end{cor}

To obtain the corollary consider $\Phi$ as a limit of
$\Phi _{n,k}$ as $k\rightarrow \infty$, see Remarks~\ref{limit}.\\ 

{\em Proof of Theorem~\ref{Utm}:}\\ 

We prove the theorem by induction on $n$. For $n=1$ the mapping $\Phi _{n,k}$ is the identity mapping 
and there is nothing to prove. Further we use the same reasoning as the one 
used in \cite{Ko4} (for $k=1$ the theorem coincides with part (2) of Theorem~1.4 in \cite{Ko4}). 
Set $P:=xQ+R$, where $R:=c_{n}(x+1)^k$. For the polynomial $xQ$ one of the numbers $a_i$ 
defined in Proposition~\ref{Kprop} equals $0$. Set 

$$ (x+1)^kxQ:=(x+1)^{n+k-1}x*(x+1)^{n+k-1}(x+h_2)*\cdots *(x+1)^{n+k-1}(x+h_{n-1})~.$$
Apply formulae (\ref{differ}). The right-hand side of the last equality is representable as 

$$ x((x+1)^{n+k-2}(x+g_2)*\cdots *(x+1)^{n+k-2}(x+g_{n-1}))~,~~{\rm where}~~g_i=\frac{(n+k-1)h_i+1}{n+k}~.$$
The last composition (excluding the factor $x$) is the representation of the polynomial $(x+1)^kQ$
in the form (\ref{Kformula}). By inductive assumption, if $\sigma _j^0$ (resp. $\sigma _j^1$ or $\sigma _j^2$)
stands for the $j$th elementary symmetric polynomial of the quantities $g_i$ 
(resp. $l_j:=(n+k)g_i/(n+k-1)$ or $h_i$), then $(-1)^j\sigma _j^0\geq 0$ (resp. $(-1)^j\sigma _j^1\geq 0$). 
We set $\sigma _0^0=\sigma _0^1=\sigma _0^2=1$. Having in mind that $h_i=l_i-1/(n+k-1)$ 
and that the signs of $\sigma _j^1$ alternate, one sees that $(-1)^j\sigma _j^2\geq 0$. Indeed, one has 
$\sigma _j^2=\sum _{\nu =0}^j(-1)^{\nu}r_{\nu}\sigma _{\nu}^1$ for some positive constants $r_{\nu}$. 
Hence $\Phi _{n,k}$ maps $U_{n}\cap \{ c_{n}=0\}$ into itself.

We show for the half-axis $Oc_{n}$ (positive for odd and negative for even $n$) that $\Phi _{n,k}(Oc_n)\subset U_n$. 
As $\Phi _{n,k}$ is affine, this implies 
$\Phi (U_{n})\subset U_{n}$.

The first $n$ coefficients of $R$ are $0$, therefore $n$ of the numbers $a_i$ defined for $\Phi _{n,k}[R]$ 
equal $\infty$ and $-s/(n+k-s)$, $s=n+k-1,\ldots ,k+1$, see Lemma~\ref{zero}. By Proposition~\ref{mult}, the remaining $k-1$ of them equal $1$. 
Therefore the numbers $a_i$ define a polynomial of the form 
$(x+1)^k(c_1^0x^{n-1}+\cdots +c_{n}^0)$ with $(-1)^{\nu}c_{\nu}^0>0$.~~~~~$\Box$

\begin{rem}
{\rm One can deduce from the proof of Theorem~\ref{Utm} that if $A\in \partial U_{n}$ (the boundary of $U_{n}$), 
then $\Phi _{n,k}(A)\in \partial U_{n}$ if and only if $A\in \{ c_{n}=0\}$.}
\end{rem}

\begin{theorem}\label{Pitm}
If $P$ is real and with $\nu$ positive roots, then at least $\nu$ of the numbers $a_i$ 
defined by formula (\ref{Kformula}) are 
negative and belonging to different intervals of the kind $I_{n,k,s}:=[-(s+1)/(n+k-1-s),-s/(n+k-s)]$. 
\end{theorem}

{\em Proof:}\\ 

The polynomial $P$ has $\nu$ positive roots. By the Descartes rule, there are at least $\nu$ sign 
changes in the sequence $\tilde{\Sigma}$ of its coefficients. On the other hand, 
when the polynomial $K_{n,k;a_i}$ is real (i.e. when $a_i$ is real), there is at most 
one sign change in the sequence of its coefficients. This follows from formula 
(\ref{Keq}) -- the numbers $((n+k-s)/(n+k))a_i+(s/(n+k))$ for $s=0,\ldots ,n+k$ form 
an arithmetic progression. For a couple of polynomials $K_{n,k;a_i}$, $K_{n,k;\bar{a}_i}$  
their SSC is a polynomial with all coefficients positive. The same is true for couples of 
polynomials $K_{n,k;a_i}$, $K_{n,k;a_j}$ with $a_i$ and $a_j$ belonging to one and the same interval 
$I_{n,k,s}$, and for polynomials $K_{n,k;a_i}$ with $a_i>0$. Hence the $\nu$ sign changes in 
the sequence $\tilde{\Sigma}$ are due only to numbers $a_i$ belonging to different intervals 
$I_{n,k,s}$.~~~~~$\Box$

\begin{rems}\label{posneg}
{\rm When $P$ or $P_m$ is hyperbolic (i.e. with all roots real), the mapping $\Phi _{n,k}$ (resp. $\Phi$) exhibits 
different properties in the cases when all roots are positive and when they are all negative. 
For instance, if all quantities $a_i$ are positive, then the composition  
$K_{n,k;a_1}*\cdots *K_{n,k;a_{n}}$ is a polynomial with all roots 
negative; this follows from Proposition~1.5 in \cite{KoSh}. But it is not true 
that when $P$ has all roots negative, then all quantities $a_i$ are real positive. 
Example: 

$$ (x+1)^{k+1}x*(x+1)^{k+1}x=(x+1)^kx(x+1/(k+2)) ~.$$
Perturb the composition factors in the left-hand side into $(x+1)^{k+1}(x\pm \varepsilon i)$. The 
polynomial to the right will have all roots negative (one of which by Proposition~\ref{mult} 
is a $k$-fold root at $-1$). This follows from the comparison of the signs 
of the constant terms to the left and right.
A similar example can be given about the mapping $\Phi$:

$$ e^x(x+1)*e^x(x+1)=e^x(x^2+3x+1)~.$$
Here $x^2+3x+1$ has two negative roots. After this perturb the two composition factors to 
the left into $e^x(x+1\pm \varepsilon i)$. For $\varepsilon >0$ small enough the polynomial 
multiplying $e^x$ in the right-hand side still has two negative roots. 

When all roots of $P$ are positive, then all quantities $a_i$ are negative, see Theorem~\ref{Pitm}. 
But when all quantities $a_i$ are negative, then all roots of $P$ are not necessarily positive. 
E.g. the following polynomial has two complex conjugate roots: 

$$(x+1)^{k}(x^2-(2kx)/(k+2)x+1)=(x+1)^{k+1}(x-1)*(x+1)^{k+1}(x-1) ~.$$
In the case of the mapping $\Phi$ an analogous example is given by the equality}

$$ e^x(x-1)*e^x(x-1)=e^x(x^2-x+1)$$
{\rm and the analog of Theorem~\ref{Pitm} in the case of the mapping $\Phi$ 
is Corollary~\ref{Descartescor} below.}
\end{rems}

\begin{nota}\label{Xinota}
{\rm For a polynomial $P=x^{n}+c_1x^{n-1}+\cdots +c_{n}$ we set}
$$ \Xi [P]:=x(x-1)\cdots (x-n+1)+c_1x(x-1)\cdots (x-n+2)+\cdots +c_{n-1}x+c_{n}~.$$
\end{nota}

\begin{rem}\label{Xirem}
{\rm One checks directly that $e^xP(x)=\sum _{j=0}^{\infty}\Xi [P](j)x^j/j!$~. It is easy to show  
that the numbers $-a_j$ 
defined by (\ref{composePhi}) are roots of the polynomial $\Xi [P]$.}
\end{rem}

Set $\Xi ^{\nu}[P]:=x^{n}+c_{1,\nu}x^{n-1}+\cdots +c_{n-1,\nu}x+c_{n,\nu}$, $c_{0,\nu}:=1$. 
It is clear that $c_{n,\nu}=c_{n,0}$ for all $\nu$.

\begin{prop}\label{iter}
1) For each real polynomial $P$ as in Notation~\ref{Xinota} there exists $\nu _0\in \mathbb{N}$ such that 
for $\nu \geq \nu _0$ the signs of $c_{0,\nu}$, $c_{1,\nu}$, $\ldots$, $c_{n-1,\nu}$ alternate. 

2) One has $\lim _{\nu \rightarrow \infty}|c_{s,\nu}/c_{s-1,\nu}|=\infty$ for $s=1$, $\ldots$, $n-1$. 

3) For $\nu$ large enough the signs of the first $n$ coefficients of $T[e^xP]$ alternate.
\end{prop}

{\em Proof:}\\ 

Observe first that $c_{0,\nu}=1$ and $c_{n,\nu}=c_{n,0}$ for all $\nu$. 
The coefficient $c_{1,\nu}$ equals $c_1-\nu n(n-1)/2$. Hence 
for $\nu$ sufficiently large this coefficient is $<0$. Moreover, after its sign stabilizes, 
its absolute value increases with each new iteration of $\Xi$ and tends to $\infty$. Hence 
$\lim _{\nu \rightarrow \infty}|c_{1,\nu}/c_{0,\nu}|=\infty$.

Suppose that each of the coefficients 
$c_{j,\nu}$, $j=1$, $\ldots$, $l-1$ of $\Xi ^{\nu}[P]$ has the {\em Property~A}: 
{\em For $\nu$ large enough its sign is the same as the one of $(-1)^j$; moreover, after its sign stabilizes, 
its absolute value increases with each new iteration of $\Xi$.}

Set $x(x-1)\cdots (x-n+1+l):=\sum _{j=0}^{n-l-1}r_{j,l}x^{n-l-j}$. Hence $r_{j,l}=0$ for 
$j>n-l-1$ and $(-1)^{j}r_{j,l}>0$. 
In particular, $r_{0,l}=1$. The constants $r_{j,l}$ depend on $n$, $l$ and $j$, but not on $\nu$. One has 
$c_{l,\nu +1}=c_{l,\nu}+\sum _{j=1}^{l-1}r_{j,l}c_{l-j,\nu}~(*).$
For $\nu$ sufficiently large the signs of all summands to the right are the same. 
Hence the coefficient $c_{l,\nu}$ also has the 
Property~A if $l<n$. This implies part 1).

Notice that $|r_{j,l}|\geq 1$ with equality only for $j=0$ and for $l=n-2$. Therefore 
part 2) of the proposition follows from (*).
Part 3) results from part 2).~~~~~$\Box$

\begin{prop}\label{Xiprop}
If the real polynomial $P$ is with all roots real positive, then the polynomial 
$\Xi [P]$ is with all roots real positive and distinct.
\end{prop}

{\em Proof:}\\ 

The non-degenerate affine mapping $\Phi$ is the limit as $k\rightarrow \infty$ of the 
non-degenerate affine mappings $\Phi _{n,k}$, see 
Remark~\ref{limit}. For each $(n,k)$ fixed the numbers $a_i$ defined for the polynomial $P$ 
from the composition product (\ref{Kformula}) are negative, see Theorem~\ref{Pitm}. 
Therefore their limits are nonpositive. The limits are $\neq 0$, otherwise 
one should have $P(0)=0$. By Remark~\ref{Xirem} the roots 
of $\Xi (P)$ are all positive. 

Set $\kappa _{a_1}*\cdots *\kappa _{a_j}:=e^xP_j(x)$. Hence $e^xP_{m-1}(x)*e^x(1+x/a_m)=e^xP_m(x)$ and  

\begin{equation}\label{sgn}
P_m(x)=(1+x/a_m) P_{m-1}(x)+(x/a_m)P_{m-1}'(x)~.
\end{equation}
By inductive assumption the polynomial $P_{m-1}$ is with distinct positive roots. Hence the 
term $(x/a_m)P_{m-1}'(x)$ changes sign at the consecutive roots of $P_{m-1}$; that is, there is a 
root of $P_m$ between any two consecutive roots of $P_{m-1}$. This makes $m-2$ distinct positive 
roots of $P_m$. One has 
sgn$(P_j(\infty ))=(-1)^j$, $j=m-1, m$ (because the quantities $a_j$ are negative) 
and sgn$P_j(0)=1$ (see (\ref{sgn})). This means that there is a root of $P_m$ in 
$(0,\lambda )$ and there is a root in $(\gamma ,\infty )$, where $\lambda$ is the smallest and 
$\gamma$ is the largest of the roots of $P_{m-1}$. Thus $P_m$ has $m$ 
distinct positive roots.~~~~~$\Box$\\ 

The following theorem (proved at the end of the paper) 
extends the Descartes rule to functions which are products 
of exponential functions and polynomials. 

\begin{theorem}\label{Descartes}
If the real degree $m$ polynomial $P$ has $k$ positive roots, $1\leq k\leq m$, then there are 
at least $k$ sign changes in the sequence of the coefficients of $T[e^xP]$. 
\end{theorem}

\begin{cor}\label{Descartescor}
If there are $k$ sign changes in the sequence of coefficients of $T[e^xP]$, then 
at least $k$ of the numbers $a_i$ 
in the composition formula (\ref{composePhi}) are negative, distinct and belonging to 
different intervals of the kind $[-l-1,-l]$, $l\in \mathbb{N}\cup \{ 0\}$. 
The conclusion is true in particular when the real polynomial $P$ has $k$ positive roots.
\end{cor}

{\em Proof:}\\ 

Formula (\ref{kappa}) implies that there is at most one change of sign in the sequence of 
coefficients of the Taylor series $T[\kappa _{a_j}]$. This change occurs only if $a_j<0$. By 
Theorem~\ref{Descartes} there are at least $k$ sign changes in the sequence of coefficients 
of $T[e^xP_m]$. Composition factors $\kappa _{a_j}$ with complex $a_j$ 
are present in (\ref{composePhi}) only in complex conjugate couples (see part 3) of Remarks~\ref{limit}), 
and for each composition of the kind $\kappa _{a}*\kappa _{\bar{a}}$ all coefficients of 
$T[\kappa _{a}*\kappa _{\bar{a}}]$ are positive. The same is true for two composition factors 
whose numbers $a_{j_1}$, $a_{j_2}$ belong to one and the same interval $[-l-1,-l]$. 
Hence the sign changes can come only from composition factors 
with negative numbers $a_j$ which belong to different intervals $[-l-1,-l]$.~~~~~$\Box$\\

\begin{cor}\label{limithyp}
For $P$ as in Notation~\ref{Xinota} there exists $\nu _0\in \mathbb{N}$ such that for $\nu \geq \nu _0$ 
the polynomial $\Xi ^{\nu}[P]$ is with real and distinct roots, $n-1$ or all of them being positive.
\end{cor}

The corollary follows from part 3) of Proposition~\ref{iter} and from Corollary~\ref{Descartescor}. 
Whether all roots or all but one are positive depends on the sign of the constant term of the polynomial. 
Indeed, the mapping $P\mapsto \Xi [P]$ preserves the constant term.

\begin{rem}
{\rm By analogy with the proof of part (5) of Theorem~1.4 in \cite{Ko4} one can prove that for each 
$(n,k)$ fixed there exists $\nu (n,k)$ such that for $\nu _0\geq \nu (n,k)$ the mapping $\Phi _{n,k}^{\nu _0}$ sends 
each point of $U_{n}$ into $U_{n}\cap \Pi _{n}$. In \cite{Ko4} this is proved for $k=1$. The following example shows 
that this is not true for the mapping $\Phi$.}
\end{rem}

\begin{ex}
{\rm Represent $f:=e^x(1+ax+bx^2)$ ($a<0$, $b>0$) in the form 
$e^x(1+x/\alpha )*e^x(1+x/\beta )$. 
Then $1/\alpha +1/\beta =a-b$, $1/(\alpha \beta )=b$. Hence 
$\Phi [f]=e^x(1+(a-b)x+bx^2)$ and  
$\Phi ^s[f]=e^x(1+(a-sb)x+bx^2)$. For every $s_0\in \mathbb{N}$ one can find 
$a<0$ and $b>0$ such that $1+(a-sb)x+bx^2$ is hyperbolic for $s\geq s_0$ 
and not hyperbolic for $s<s_0$.}
\end{ex}

\begin{prop}\label{notV}
For $m=3$ the mapping $\Phi$ does not send the set $V_{m}$ into itself.
\end{prop}

{\em Proof:}\\ 

Consider the functions of the kind $e^x(x^3+ax^2+bx+c)$, $a\leq 0$, $b\geq 0$, $c\leq 0$. 
Their subset whose roots have non-negative real parts is bounded by the hyperbolic paraboloid 
$\mathcal{P}:c=ab$ and the hyperplanes $\mathcal{H}_1:a=0$, $\mathcal{H}_2:b=0$ and $\mathcal{H}_3:c=0$. 
It is defined by the system $c\geq ab$, $a\leq 0$, $b\geq 0$, $c\leq 0$. Its boundary is 

$$ A\cup B~~,~~{\rm where}~~A:=\{ ~c=0~,~a\leq 0~,~b\geq 0~\} ~~,~~B:=\{ ~c=ab~,~a\leq 0~,~b\geq 0~\} ~.$$
The polynomials corresponding to the set $A$ have a root at $0$, the ones from $B$ are of the form 
$R:=(x-d)(x^2+\Lambda )=x^3-dx^2+\Lambda x-d\Lambda $, $d\geq 0$, $\Lambda \geq 0$. Set 

\begin{equation}\label{abc}
e^x(x^3-dx^2+\Lambda x-d\Lambda )=e^x(x+\alpha )*e^x(x+\beta )*e^x(x+\gamma )~,
\end{equation}
$\sigma _1:=\alpha +\beta +\gamma$, $\sigma _2:=\alpha \beta +\alpha \gamma +\beta \gamma$, $\sigma _3:=\alpha \beta \gamma$. 
Comparing the coefficients of $1$, $x$ and $x^2$ in the two sides of (\ref{abc}) one obtains the system 

$$ \sigma _3=-d\Lambda ~,~1+\sigma _1+\sigma _2+\sigma _3=\Lambda
-d\Lambda ~,~8+4\sigma _1+2\sigma _2+\sigma _3=-2d+2\Lambda -d\Lambda
~.$$
This means that $\Phi [e^xR]=e^x(x^3+(-d-3)x^2+(\Lambda +d+2)x-d\Lambda )$. The coefficients $a$, $b$, $c$ of 
the last polynomial factor satisfy the condition $c=(a+3)(b+a+1)$. This defines a hypersurface $Y\subset \mathbb{R}^3$. 
Consider the intersection $B\cap Y$. It is defined by the conditions $a\leq 0$, $b\geq 0$, $c=ab=(a+3)(b+a+1)$. 
The point $W:=(a,b,c)=(-2,1/3,-2/3)$ belongs to this intersection. Fix $a=-2$ and vary $b$. Close to the point $W$ 
there are points of $Y$ which are inside and points which are outside the domain $V_3$. This proves the proposition.~~~~~$\Box$\\

{\em Proof of Theorem~\ref{Descartes}:}\\ 

$1^0$. Assume (which is not restrictive) that $P$ is monic and that $P(0)\neq 0$. 
Set $P:=x^m+d_1x^{m-1}+\cdots +d_m$, $d:=|d_1|+\cdots +|d_m|$. 
There exists $N\in \mathbb{N}$ such that the coefficient of $x^j$ in $T[e^xP]$ 
is positive for $j\geq N$. Indeed, this coefficient equals 

$$ 1/(j-m)!+d_1/(j-m+1)!+\cdots +d_m/j!$$
which is $>0$ for $d<j-m+1$ (i.e. for $j>d+m-1$). 

$2^0$. Suppose that $P$ has a  
root $x_0>0$ of multiplicity $\mu >1$ (if all positive roots are simple, then go directly to $3^0$). 
Denote by $x_1$, $\ldots$, $x_s$ its other positive roots and by $m_1$, $\ldots$, $m_s$ their multiplicities. 
For $\varepsilon >0$ small enough the coefficients of $T[e^xP+\varepsilon x^N(\prod _{i=1}^s(x-x_i)^{m_i})(x-x_0)^{\mu -1}]$ 
have the same signs as the coefficients of $T[e^xP]$. The root $x_0$ bifurcates into 
a root of multiplicity $\mu -1$ and a simple root close to it, both positive. The other positive roots and their 
multiplicities remain the same. In the same way one can change the function $e^xP$ to a nearby one 
$e^xP+g$ ($g$ is a polynomial) with 
the same number of positive roots (counted with multiplicity, but which are all simple) and the same signs of 
its Taylor coefficients.

$3^0$. Fix an interval $I:=[\delta _1,\delta _2]$ ($0<\delta _1<\delta _2$) 
containing in its interior all positive roots of $e^xP+g$. 
The series $T[e^xP]$ converges absolutely for 
all real $x$ and its coefficients except finitely many are positive. Therefore one can find a partial sum $S$ of $T[e^xP+g]$ with the same number of roots in $I$ as $e^xP+g$ (all of them being simple) and with the same number of sign changes 
in the sequence of its coefficients. Hence $S$ has $\geq k$ positive roots and the number of sign changes 
in the sequence of the Taylor coefficients (which is the same for $T[e^xP]$) is $\geq k$.~~~~~$\Box$

\end{document}